\documentclass{amsart}
\usepackage{amsmath}
\usepackage{amsfonts}
\usepackage{color}
\usepackage{amssymb,mathrsfs}
\usepackage{t1enc}
\usepackage[latin1]{inputenc}
\usepackage[english]{babel}
\usepackage[numbers]{natbib}
\usepackage{amssymb,stmaryrd}
\usepackage{caption}
\usepackage{url}

\usepackage[colorlinks=false, linktocpage=true,backref,pagebackref,hidelinks]{hyperref}
\DeclareMathSymbol{\varGamma}{\mathord}{letters}{"00}
\DeclareMathSymbol{\varPi}{\mathord}{letters}{"05}
\DeclareMathSymbol{\varLambda}{\mathord}{letters}{"03}
\newcommand{\Cs}{\mathscr{C}}
\newcommand{\Ps}{\mathscr{P}}

\newcommand{\beq}{\begin{equation}}
\newcommand{\eeq}{\end{equation}}
\newcommand{\beqarr}{\begin{eqnarray}}
\newcommand{\eeqarr}{\end{eqnarray}}
\newcommand{\beqa}{\begin{eqnarray*}}
\newcommand{\eeqa}{\end{eqnarray*}}
\newtheorem{theorem}{Theorem}

\newtheorem{definition}{Definition}
\newtheorem{remark}{Remark}

\begin{document}
\thispagestyle{empty}
\renewcommand{\thefootnote}{\fnsymbol{footnote}}

\title{ { Light cone and Weyl }   compatibility of conformal and projective structures  }
\author{Vladimir S. Matveev,  Erhard Scholz}

\address{V. \ M.:  Mathematics Institute, Faculty of Mathematics and Informatics,
Friedrich-Schiller-Universit\"at Jena, 07737 Jena, Germany
\\ E.\ S.: 
Faculty of  Math./Natural Sciences, and Interdisciplinary Centre for Hist. and Phil. of Science, University of Wuppertal, Germany}

\email{vladimir.matveev@uni-jena.de, scholz@math.uni-wuppertal.de}

\renewcommand{\thefootnote}{\arabic{footnote}}

%%\author{\normalsize Vadimir Matveev,\footnote{Friedrich-Schiller-Universit\"at Jena, Faculty of Mathematics,   Germany,  vladimir.matveev@uni-jena.de} \quad Erhard Scholz\footnote{University of Wuppertal, Faculty of  Math./Natural Sciences, and Interdisciplinary Centre for Hist. and Phil. of Science, Germany, \;  scholz@math.uni-wuppertal.de}}
\date{12/25/2019 } 

%\vspace{2.5mm}
\begin{abstract}
In the literature different concepts  of compatibility between a projective structure $\Ps$ and a conformal structure $\Cs$ on a differentiable manifold are used. In particular compatibility in the sense of Weyl geometry is slightly more general than compatibility in the Riemannian sense. In an often cited  paper  \citep{EPS}  Ehlers/Pirani/Schild introduce still another criterion which is natural from the physical point of view: 
 every { light like   geodesics of $\Cs$ is  a    geodesics} of $\Ps$.  Their claim that this type of compatibility is sufficient for introducing a Weylian metric  has recently been questioned \citep{Trautman:EPS,Matveev/Trautman,Scholz:2019}. Here it is proved that the conjecture of EPS is correct. 
\end{abstract}
\maketitle

%%%%%%%%%%%%%%%%%%%%%%%%%%%
\section*{\small Introduction}
\addcontentsline{toc}{section}{\protect\numberline{}Introduction}
In a widely read paper \citep{EPS} J. Ehlers, F. Pirani and A. Schild (EPS)   argued  that a  projective structure $\Ps$ and a conformal structure $\Cs$ on a differentiable manifold $M$ determine a Weylian metric on $M$, if only the  { geodesics  } of $\Ps$ are 
 {light like  geodesics}  of  $\Cs$. This statement complements a proposal of H. Weyl for basing the geometrical framework of   gravity theory on the observable structures of particle trajectories and light propagation, rather than on the chronogeometric behaviour of clocks or  rods  \citep{Weyl:projektiv_konform}. EPS claimed that the above mentioned  light-cone condition for $\Ps$ and $\Cs$  is sufficient for securing  the existence of a Weylian metric  which Weyl had  assumed from the outset; but the argument given to substantiate the statement remains vague and is far from a mathematical proof \citep{Trautman:EPS}. 
 The aim of the present paper is to fill in the gap and to prove the central statement of EPS.

  The paper is structured as follows. We start with  short remarks on the history of the problem (sec. \ref{section Hist}). After clarifying the central  concepts  involved in the question (sec. \ref{section Definitions})  { we   prove   the EPS conjecture }(sec. \ref{section EPS conjecture}). Finally we discuss why  we think it matters (sec. \ref{section Discussion}).

\vspace{1ex}
\noindent {\bf Acknowledgements.}  We thank Sergei Agafonov, Thomas Mettler and
Andrzej  Trautman for useful comments. The work was started  during the conference ``New methods in differential geometry'' supported   by the DAAD (via Ostpartnerschaft programme) and by the University of Jena. 
V.\, Matveev  thanks the DFG for partial support via projects MA 2565/4
and MA 2565/6.

\vspace{1ex} 
\section{\small Some historical remarks \label{section Hist}}
In 1918 H. Weyl generalized the concept of a Riemannian metric  in order to avoid the possibility of direct metrical comparison of vectors or other fields at finitely distant points \citep{Weyl:InfGeo,Weyl:GuE}. He   introduced  a scale gauge connection in addition to a conformal structure,  thus defining what later would be  called a {\em Weylian metric} on a differentiable manifold (see below, Def. \ref{def Weyl structure}(ii)). 
In his book {\em Space - Time - Matter} and lectures given at Barcelona in 1922 he argued that the geometry  of special relativity, i.e. the affine structure of Minkowski space, can be established without rods and clocks by projective geometry and the specification of a quadric at infinity in the sense of Klein's Erlangen program \citep{Weyl:STM,Weyl:ARP1923}. Generalizing this idea  he  argued  that the geometrical structure of the general theory of relativity (GTR) can be based on the  mathematical description of the inertial motion of test particles and light rays, rather than on the behaviour of rods and clock. He was able to underpin this view by showing that the generalized metrical structure of gravity, which  he had proposed, i.e. a Weylian metric on $M$,  is uniquely determined by its  associated projective and and conformal structures \citep{Weyl:projektiv_konform}.\footnote{See also Weyl's discussion with Einstein in 1918  \citep[vol. 8B]{Einstein:CP}. } 
He did not discuss, however, the conditions under which a projective structure and a conformal structure determine a Weylian metric.

About half a century later, in 1972, 
 J. Ehlers, F. Pirani and A. Schild, sketched an even more ambitious program for establishing the fundamental conceptual  framework  of general relativity. They wanted to base  even the    differentiable structure on the set of  spacetime events on more general, physically more or less plausible, assumptions (called ``axioms'') on the relation between events, particle paths and light propagation in spacetime \citep{EPS}. 
  This was considered as an attempt for a physically motivated ``constructive axiomatics'' of GTR and included  the central claim that  a projective structure and a conformal structure which are compatible on the light cones  (see below Def. \ref{def compatibility}(i)) determine a Weylian metric on the spacetime manifold.\footnote{More precisely EPS speak of a ``Weyl space'' if a light cone compatible pair of projective and conformal structures is given. In the rest of the paper they suggest that a ``Weyl space'' can be endowed with a Weylian metric.}
  The EPS paper led to a series of follow up investigations which  in many  cases concentrated on conditions which would reduce the Weyl geometric structure to a Riemannian (Lorentzian) one, often introducing additional information of quantum physics (Dirac field, complex scalar field)   \citep{Audretsch:1983,AGS,Audretsch/Hehl/Laemmerzahl,Audretsch/Laemmerzahl:1994,Coleman/Korte:inertial_conformal}.\footnote{For more details see \citep{Scholz:2018Resurgence}.}
 In these investigations the arguments of EPS were usually accepted, although the authors had qualified their arguments as not necessarily mathematically satisfying.\footnote{ ``A fully rigorous formalization  has not yet been achieved, but we nevertheless hope that the
main line of reasoning will be intelligible and convincing to the sympathetic reader'' \citep[p. 69f.]{EPS}.
}

The EPS paper  was republished in 2012 as a by then classical text (``Golden Oldie'') with an editorial comment by A. Trautman \citep{Trautman:EPS}.   In his comments  the editor  raised doubts with regard to the status of the existence statement of EPS for a Weylian metric.  He made clear that the arguments given in the original paper were rather  vague and  far from a mathematical proof. The existence statement  ought thus to be considered a conjecture rather than a theorem as which it had been treated in large parts of the literature up to then.
A first investigation of the case in a joint paper of Trautman with one of the present authors draws the conclusion that the EPS statement is wrong \citep{Matveev/Trautman}. This judgment is based, however, on the criterion of Riemann compatibility between  projective and conformal structures and thus on  a too narrow understanding of the Weyl geometric setting. In the following argument it will be shown that the EPS conjecture is, in fact, true.

\vspace{1ex} 
\section{\small Definitions \label{section Definitions}}
We consider smooth  manifolds and maps of class at least \(\Cs^2\). All geometric objects
on an \(n\ge 3 \)-dimensional manifold are referred to local coordinates \((x^i)\),\; \(i=1,\dots,n\).

A \emph{conformal structure} on a manifold \(M\) is an equivalence class \(\Cs\)  of
 metric tensors \(g\) with respect to the following equivalence relation
\[
g \sim  g'\;\Longleftrightarrow \mbox{there is a function \(\varphi\) on \(M\) such that \(g'=g \exp
{2\varphi}\).}
\]
If \(g\in\Cs\), then \(\Cs\) can be denoted by \([g]\).
We assume that the metric has indefinite signature, since otherwise the compatibility condition 
(see  below) is empty.

Two  symmetric linear connections
 \(\varGamma=(\varGamma^{i}_{jk})\)  and   \(\varGamma'=({\varGamma'}^{i}_{jk})\) are said to be \emph{projectively equivalent} if their {  geodesics coincide.  Here and below we consider geodesics without preferred parameterization (in literature they are sometimes  called \emph{autoparallel curves}).}   Projective equivalence  is clearly an equivalence relation on the set of all  symmetric linear connections on \(M\). An equivalence class \(\Ps\) with respect to this relation is  called a {\em projective structure};  it is denoted by \([\varGamma]\) if it contains \(\varGamma\). It    can be formulated as the condition
 \begin{center}
\(
   \varGamma \sim \varGamma'\in\Ps\;\
   \Longleftrightarrow\; \mbox{there is a 1-form \(\psi\) so that \(\varGamma^{\prime i}_{jk}=\varGamma^{i}_{jk}+\delta^i_j\psi_k+
   \delta^i_k\psi_j.\)}
   \)
\end{center}

We consider here a question of Weyl's generalization of Riemannian geometry proposed in \citep{Weyl:InfGeo,Weyl:GuE,Weyl:STM}. In the more recent literature this type of generalization has been formulated for various differential geometric structures \citep{Higa:1993, Ornea:2001,Gilkey_ea},\footnote{For a concept of Weyl structures in the context of Cartan geometries modeled after a pair $(G;P)$ with $P$ a parbolic subroup of the Lie group $G$ see \citep[chap. 5]{Cap/Slovak:2009}.}
 We use it in the sense of semi-Riemannian Weyl structures (Def. \ref{def Weyl structure} (i)) which are close to the Weylian manifolds (Def. \ref{def Weyl structure}  (ii)) considered by Weyl himself. 
\begin{definition} \label{def Weyl structure}
\begin{itemize}
\item[(i)]
 A  (semi-Riemannian) {\sc Weyl structure} is given by triple $(M, \Cs, \nabla)$
where   $M$ is  a differentiable manifold, $\Cs =[g]$  a {\em conformal class}  of (semi-) Riemannian metrics $g$ on $M$,  and $\nabla = \nabla(\Gamma)$ the {\em covariant derivative}  of  a  torsion free affine connection $\Gamma$ , constrained by the {\em compatibility condition} that for any $g\in \Cs$ there is a differential 1-form $\varphi_g$ s.th. 
$ \nabla g + 2 \varphi_g \otimes g  = 0 $.
\item[(ii)]  A {\sc Weylian manifold}  $(M,[(g,\varphi)]  )$  is  a differentiable manifold $M$ endowed with a {\sc Weylian metric}  defined by  an equivalence class of pairs $(g,\varphi)$, where $g$ is a (semi-) Riemannian metric  and $\varphi$  a (real valued) differential 1-form on $M$. Equivalence is defined by conformal rescaling   $g \mapsto \tilde{g} =\Omega^2 g$ and the corresponding gauge transformation for $\varphi \mapsto \widetilde{\varphi}=\varphi- d \ln \Omega$.
\end{itemize}
\end{definition}
Weyl showed that any Weylian metric has a uniquely determined compatible affine (i.e. symmetric linear) connection $\Gamma(g,\varphi)$, 
\beq  \varGamma(g,\varphi)^i_{jk} = \digamma^i_{jk} +\delta^i_j\varphi_k + \delta^i_k \varphi_j -g_{jk}\varphi^i \, ,\label{Gamma-Weyl}
\eeq 
where $ \digamma$ denotes the Levi-Civita connection of $g$. It is 
 independent of the  representative $(g, \varphi)$ of the Weylian metric \citep{Weyl:InfGeo}. {\em Metric compatibility} in the sense of Weyl geometry means that the  the lengths of vectors parallel transported by $\varGamma(g,\varphi)$  and measured in $g$ change infinitesimally with $\varphi$. In streamlined form this  means that for the covariant derivative $\nabla = \nabla(g,\varphi)$ defined by  $\Gamma(g,\varphi)$ the following holds:
\beq \nabla g + 2 \varphi \, g = 0 
\eeq 
Taking into account the gauge transformation for $\varphi$  the definitions (i) and (ii) above turn out to be  equivalent.
 
The compatibility of a projective structure with a conformal structure can now be considered from different perspectives. We use the following terminology:
\begin{definition} \label{def compatibility}
We say that a projective structure \(\Ps=[\varGamma]\) and  a conformal structure  \(\Cs=[g]\) are 
\begin{itemize}
 \item[(i) ]{\sc light cone compatible} 
 if any light-like  geodesic of $g\in \Cs$ is an auto-parallel for some $\varGamma \in \Ps$;
 \item[(ii)] {\sc Riemann compatible} if there is $ g \in\Cs$ such that its Levi-Civita connection  $\digamma(g) \in  \Ps$; 
 \item[(iii)] {\sc Weyl compatible}  if for some $g\in \Cs$  a differential 1-form $\varphi$ can be found such that the  affine connection $\varGamma (g,\varphi)$ of the Weylian metric $[(g,\varphi)]$ satisfies $\varGamma (g,\varphi) \in \Ps$.\footnote{If  this holds for  some $g\in \Cs$, then for {\em any} $g\in \Cs$.} 
 \end{itemize}
\end{definition}

 \begin{remark}  (i)   is  independent of the choice of the connection from the projective class of \(\varGamma\) and of the choice of the metric from the conformal class of $g$; it is used in \citep{EPS}.  We also use the abbreviation {\em compatibility} without further specification for (i).   In the context of Riemannian geometry (ii) appears most natural; it is used  also in \citep{Matveev/Trautman}. (iii) is a straight forward generalization of Riemann compatibility to the context of Weyl geometry and was implicitly considered by Weyl in  \citep{Weyl:projektiv_konform}. 
 \end{remark}

Weyl compatibility implies (light cone) compatibility  \citep{Weyl:InfGeo},  similarly so for Riemann compatibility. On the other hand,  light cone compatibility does not imply  Riemann compatibility.  The question remains whether light cone compatibility is strong enough to imply Weyl compatibility
{ (assuming that the metric has indefinite signature so the light cone exists).}

A trivial example of (light cone) compatible projective and conformal structures is as follows: take 
any two $1-$forms $\varphi= \varphi_i$ and $\eta= \eta_i$, denote by $\digamma^i_{jk}$ the Levi-Civita connection of any $g \in \Cs$ and consider $\Ps=[\varGamma ]$ with
\begin{equation} \label{new}
\varGamma^{i}_{jk} = \digamma^i_{jk} + \varphi^i g_{jk} +  \eta_j \delta^i_k +  \eta_k \delta^i_j.
\end{equation}
$\varGamma$ and $\digamma$ are obviously  projectively equivalent on the light cones of $\Cs$; and this property does not depend on the choice of the representative of $\Ps$.

{
In our paper we prove that any pair  $(\Ps, \Cs)$ of 
 compatible  projective and  conformal structures are related by the formula \eqref{new}, see Theorem \ref{below}. This implies that  
 compatible $\Ps$ and $\Cs$ are  also  Weyl compatible, because  (\ref{Gamma-Weyl}) shows that in this case also the invariant affine connection of the Weylian metric $[(g,\varphi)]$ is  projectively equivalent to $\varGamma$. }

\vspace{1ex} 
\section{\small  Proof of the EPS conjecture  \label{section EPS conjecture}}

\begin{theorem} \label{below} Let $g$ be a metric of indefinite signature  on $\mathbb{R}^n$ with $n\ge 3$. 
   If $[\varGamma^{i}_{jk}]$ is  compatible with $[g]$, then \eqref{new} holds. 
\end{theorem}

{\bf Proof.}
For any light-like geodesic $\gamma$ we have 
$$
\nabla^g_{\dot \gamma}  \dot  \gamma= 0 \ \textrm{and} \ \nabla^{\varGamma}_{\dot \gamma}   \dot \gamma= \beta(\gamma, \dot \gamma)\dot \gamma   
$$
(for some function $\beta$). 
Subtracting one equation from the other, we obtain 
\begin{equation} \label{system}
 D^i_{jk}   \dot \gamma^j \dot \gamma^k= \beta(\gamma, \dot \gamma) \dot \gamma^i , 
\end{equation}
where the ``difference''  $D$ is given by  $D^s_{jk}=  \varGamma^s_{jk}-  \digamma^s_{jk}$; it is a tensor. 

We view \eqref{system} as a system of linear equations on the components of
$D$ (assumed to be symmetric in the lower indices);  the system contains infinitely many equations since  as $\dot \gamma$ we can take any light-like vector. Our goal is to show that the general solution of this system is the one coming from \eqref{new}.

 For any vector $v=v^i$  we consider 
\begin{equation}\label{prod}
D^i_{jk} v^jv^k v^s - D^s_{jk} v^jv^k v^i. 
\end{equation}
This polynomial in $v$ of  degree 3 expression vanishes for any $v$ such that $g(v,v)=0$; since the set 
$v\in \mathbb{R}^n $ such that $g(v,v)=0$ is an irreducible quadric,  the expression  \eqref{prod} is divisible by $g(v,v)$, so it is equal to 
\begin{equation}\label{prod1}
 g(v,v)  {\omega^{is}}_p v^p
\end{equation}
for some ${\omega^{is}}_p$ skew-symmetric with respect to $i,s$. On the other side,  \eqref{prod} has the
following  property: for any two 1-forms  $\sigma_i$, $\zeta_i$  such that $\sigma_i v^i=\zeta_i v^i=0$, 
 if we contract \eqref{prod} with 
   $\sigma_i \zeta_s$ we obtain zero. Thus, at any $v$ and  for any such $\sigma$ and $\zeta$  we have  
$$
 \sigma_i \zeta_s {\omega^{is}}_p v^p=0. 
$$

Then, for any point $v$ the  contravariant 2-form  ${\omega^{is}}_p v^p$ has rank two so it 
 is given by $\varphi^k v^i - v^k \varphi^i$ for some $\varphi$. 
Let us show that this implies that \begin{equation} \label{sol} D^i_{jk} = \varphi^ig_{jk} + \delta^i_j \eta_k +  \delta^i_k \eta_j \end{equation} 
 as we want. 

In order to do this, we consider the equation 

\begin{equation}\label{prod2}
D^i_{jk} v^jv^k v^s - D^s_{jk} v^jv^k v^i =(\varphi^i v^s - v^i \varphi^s) g_{jk} v^j v^k. 
\end{equation}
We view it as an  equation on $D$ (and assume $\varphi$ is known). It should be fulfilled for all vectors  $v$. It is a
 system of linear inhomogeneous equations. The corresponding homogeneous system is $D^i_{jk} v^jv^k v^s - D^s_{jk} v^jv^k v^i=0$. It is equivalent to the condition that for  every $v$ we have that $D^i_{jk} v^jv^k$ is proportional to $v^i$. Then, its solution space  is the space  of pure trace 
 tensors (i.e., of the form $D^i_{jk}= \delta^i_{j} \eta_k + \delta^i_k \eta_j$). Now,  the general solution of an inhomogeneous system is one solution plus all solutions of the corresponding homogeneous system; since 
$D^i_{jk} = \varphi^ig_{jk}$ is a solution of the system \eqref{prod2}, the general solution is \eqref{sol} as we claimed.

\begin{remark} Arguing as  in \cite[\S 2]{Matveev/Trautman}, one can extract from $D^i_{jk} := \varGamma^{i}_{jk} - \digamma^i_{jk} = \varphi^i g_{jk} +  \eta_j \delta^i_k +  \eta_k \delta^i_j$ a formula for $\varphi^i$: indeed, 
   $$ \left(n-\frac{2}{n+1}\right)\varphi^i = \left(D^i_{jk}-\frac{1}{n+1} D^s_{sk} \delta^i_j - \frac{1}{n+1} D^s_{sj} \delta^i_k\right)g^{jk}.$$
 The Weyl structure  corresponding to this $\varphi^i$  is integrable if  $\varphi_i$ is closed. \end{remark}  
  
\begin{remark} An alternative equivalent way to solve the system \eqref{system} is as follows: 
The system is invariant with respect to  the natural action of the group $O(g)$. Then, the solution space is also invariant and so is the direct sum of irreducible subspaces. Considering all  irreducible subspaces and substituting them as ansatz in \eqref{system} (which is a standard exercise) shows that only the subspaces 
of tensors of the form  $\eta_j \delta^i_k +  \eta_k \delta^i_j$ and of  the form  $\varphi^i g_{jk}$ are solution spaces for the system. 

Our proof is more elementary and works in all dimensions and all signatures (recall  that decomposition in  irreducible subspaces w.r.t. to the action of $O(g)$ depends on the dimension). 
\end{remark}  

\vspace{1ex} 
\section{\small Discussion \label{section Discussion}}
Theorem \ref{below}  { and discussion at the end of Section \ref{section Definitions}}  show that EPS  were right in assuming that the light cone compatibility of $(\Ps, \Cs)$ is equivalent to the existence of a Weylian metric with the given projective and conformal structures.  The Weylian metric is, moreover, well determined because of Weyl's uniqeness theorem, see section \ref{section Hist}. In this sense the central claim of EPS has been  vindicated.\footnote{ We do not deal with those parts of the paper in which the differentiable structure on $M$ and the structures $\Ps,\Cs$ are derived from the rather involved ``axioms'' of \citep{EPS}.}
In principle,  a Weylian metric on spacetime can thus be read off from sufficiently detailed knowledge of the free fall trajectories of test  particles and of the gravitational bending of light.   This fact may give support for a modified gravity approach to the problem of dark matter at the astrophysical level (galaxies and galaxy clusters) \citep{Fatibene/Francaviglia:2012,Scholz:WST2019}. The central physical question  for such an approach is then whether a theoretical coherent and empirically confirmed dynamics of the underlying field content can be found.

%\subsection{\small xxx3 \label{subsection xxx3}}

%\subsection{\small xxx4 \label{subsection xxx4}}

%%%%%%%%%%%%%%%%%%%%%%%%%%%%

\end{document}